\documentclass[11pt,a4paper]{amsart}
\usepackage[top=1.15in, bottom=1.15in, left=1.17in, right=1.17in]{geometry}
\usepackage[utf8]{inputenc}
\usepackage{amssymb}
\usepackage{pdflscape}
\usepackage{amscd}
\usepackage{dsfont}
\usepackage{fancyhdr}
\input xy 
\xyoption{all}
\usepackage{amsmath,amsfonts,stmaryrd}
\usepackage[abbrev,msc-links]{amsrefs}
\usepackage{tikz}
\usetikzlibrary{shapes,arrows}
\theoremstyle{plain}

\newtheorem{thm}{Theorem}[section]
\newtheorem{lem}[thm]{Lemma}

\newtheorem{pro}[thm]{Proposition}
\newtheorem{cor}[thm]{Corollary}

\theoremstyle{definition}
\newtheorem{dfn}[thm]{Definition}
\newtheorem{exa}[thm]{Example}
\newtheorem{rem}[thm]{Remark}

\DeclareMathOperator{\Hom}{Hom}

\DeclareMathOperator{\End}{End}
\DeclareMathOperator{\Ext}{Ext}

\DeclareMathOperator{\Lamod}{\Lambda\text{-mod}}

\DeclareMathOperator{\id}{id}
\DeclareMathOperator{\idim}{id}
\DeclareMathOperator{\pdim}{pd}

\newcommand{\oTo}{\xymatrix{ \ar@{^{(}->}[r]|{\mathbf{O}}& }} % open immersion
\newcommand{\cTo}{\xymatrix{ \ar@{^{(}->}[r]|{\mathbf{|}}& }} % closed immersion
\newcommand{\coTo}{\xymatrix{ \ar@{^{(}->}[r]|{\mathbf{O}}|{\mathbf{|}}& }} %locally closed

\DeclareMathOperator{\Bild}{Im}

\newcommand{\La}{\Lambda}

\newcommand{\K}{\mathbb{K}}

\newcommand{\N}{\mathbb{N}}

\newcommand{\mcE}{\mathcal{E}}

\newcommand{\mcI}{\mathcal{I}}

\newcommand{\mcP}{\mathcal{P}}

\newcommand{\mcX}{\mathcal{X}}
\newcommand{\mcY}{\mathcal{Y}}

\DeclareMathOperator{\add}{add}

\DeclareMathOperator{\Gr}{Gr}

\usepackage{todonotes}

\begin{document}
%\title{Quasi-projective varieties are Grassmannians for fully exact subcategories of quiver representations}
\title[Quasi-projective varieties are Grassmannians for exact categories]{Quasi-projective varieties are Grassmannians for fully exact subcategories of quiver representations}

\date{\today}

\keywords{Quasi-projective variety, quiver Grassmannian, exact category}

\subjclass[2020]{18G80 14A10}

\author{Alexander P\"utz}
\address{A. P\"utz:\newline
Institute of Mathematics\\
University Paderborn\\
Warburger Str. 100\\\newline
D-33098 Paderborn\\
Germany}
\email{alexander.puetz@math.uni-paderborn.de}

\author{Julia Sauter}
\address{J. Sauter:\newline
Faculty of Mathematics\\
Bielefeld University\\
PO Box 100 131\\\newline
 D-33501 Bielefeld\\
Germany}
\email{jsauter@math.uni-bielefeld.de}

\begin{abstract} 
Reineke and independent other authors proved that every projective variety arises as a quiver Grassmannian. We prove the claim in the title by restricting Reineke's isomorphism to Grassmannians for a fully exact subcategory. 
\end{abstract}

\maketitle

%++++++++++++++++++++++++++++++++++++++++++++++++++++++++++++
\section{Introduction}
%++++++++++++++++++++++++++++++++++++++++++++++++++++++++++++
Quiver Grassmannians are varieties parameterizing subrepresentations of a fixed dimension vector for a given quiver representation. Projective varieties arise as quiver Grassmannians in several ways, as shown by Reineke \cite{R} with a precursor by Hille \cite{H} and further variants by Ringel \cites{Ri,Ri2}. %\\
%As  we introduce in Section~\ref{sec:open-subcats} open subcategories in the category of modules over basic finite-dimensional algebras (inspired by Drozd and Greul \cite{DG}).  

In Section~\ref{sec:open-subcats} we introduce open subcategories in the category of modules over a basic finite-dimensional algebra (inspired by Drozd and Greul \cite{DG}). %This is analogous to quasi-projective varieties which open subvarieties of projective varieties. 
In Corollary~\ref{openImpliesExt-closed} we show that open subcategories are always extension-closed. We consider them with the restricted exact structure as a fully exact subcategory. For example, the full 
subcategory of all modules which are right Hom-orthogonal to a module is open (Example~\ref{ex.open-subcat}). In Section~\ref{sec:exact-quiver-grass}, we introduce Grassmannians for the above open fully exact subcategories. By construction they are open subschemes in the Grassmannian of the algebra, cf. Lemma~\ref{lem:subvariety-structure}. In Section~\ref{sec:quasi-proj}, we restrict Reineke's isomorphism to prove the following:
\begin{thm} (cf. Thm~\ref{thm:quasi-projective-varieties-as-quiver-Grass}) Every quasi-projective variety is a Grassmannian for an open fully exact subcategory in the category of quiver representations.  
\end{thm}
More precisely, we show that every quasi-projective variety is isomorphic to a Grassmannian of an exact category (denoted by $\Gr_{\mcE} (\mathbf{d},V)$), where 
\begin{itemize}
\item[(*)] $Q$ is an acyclic quiver on three vertices,
\item[(*)] $\mcE$ consists of $Q$-representations right Hom-orthogonal to a brick, %(i.e. an object with endomorphism ring $\mathbb{K}$),
\item[(*)] $V$ is a brick,
\item[(*)] $\mathbf{d}$ is a thin dimension vector (i.e. $d_i\leq 1$ for all vertices $i$ of $Q$).
\end{itemize}
%A brick is a module with endomorphism ring $\mathbb{K}$ (as we assume $\mathbb{K}$ is algebraically closed). 

\subsection*{Acknowledgement}
The first author was fully and the second author partially funded by the Deutsche Forschungsgemeinschaft (DFG, German Research Foundation) -- Project-ID 491392403 -- TRR 358.  
% Ich habe das nur umgeschrieben damit es in die Zeile passt. Wenn wir uns noch bei jemandem bedanken wollen, können wir das auch wieder rückgängig machen.

%++++++++++++++++++++++++++++++++++++++++++++++++++++++++++++
\section{Open fully exact subcategories}\label{sec:open-subcats}
%++++++++++++++++++++++++++++++++++++++++++++++++++++++++++++

Let $\Lambda$ be a basic finite-dimensional algebra and let $\La$-mod %$\mathrm{mod}_K \Lambda$ 
denote the category of finite-dimensional left $\Lambda$-modules over an algebraically closed field $\K$ (see \cite[Chapter~I]{ASS06} for more detail). For an exposition on exact categories see \cite{Bu10}.

Let $Q$ be the quiver associated to $\Lambda$ with admissible ideal $I \subseteq \K Q$ so that $\Lambda \cong \K Q/I$ (see \cite[Theorem~II.3.7]{ASS06}). $Q_0$ denotes the vertices of $Q$ and $Q_1$ denotes the set of oriented edges $(\alpha : i \to j)$.
For $\mathbf{d} \in \N_{0}^{Q_0}$, $\mathrm{R}_\mathbf{d}(Q)$ is the variety of $\mathbf{d}$-dimensional $Q$-representations:
\[ \mathrm{R}_\mathbf{d}(Q) = \bigoplus_{(\alpha:i\to j) \in Q_1} \K^{d_j\times d_i}.\]
The group $\mathrm{G}_\mathbf{d} := \prod_{i \in Q_0} \mathrm{GL}_{d_i}(\K)$ acts on $\mathrm{R}_\mathbf{d}(Q)$ via conjugation: 
\[ g.M := (g_jM_\alpha g_i^{-1})_{(\alpha:i \to j) \in Q_1},\]
for $g \in \mathrm{G}_\mathbf{d}$ and $M \in \mathrm{R}_\mathbf{d}(Q)$. It contains the closed $\mathrm{G}_\mathbf{d}$-stable subvariety $\mathrm{R}_\mathbf{d}(\Lambda)$ of $\mathbf{d}$-dimensional $\Lambda$-modules. %We sometimes use $\alpha$ as shorthand notation for $(\alpha:i \to j) \in Q_1$.
\begin{dfn}\label{RdE}
For an additively closed subcategory $\mcX \subseteq \La\mathrm{-mod}$, $Q$ the quiver of $\Lambda$ and $\mathbf{d} \in \N_{0}^{Q_0}$, define
\( \mathrm{R}_\mathbf{d}(\mathcal{X}) := \{M \in \mathrm{R}_\mathbf{d}(\La) : M \in \mcX\}\). Then $\mathcal{X}$ is called 
\begin{itemize}
\item[(i)] \textbf{open} 
if for every dimension vector $\mathbf{d}\in \N_{0}^{Q_0}$, the subset $\mathrm{R}_\mathbf{d}(\mcX)$ is open in $\mathrm{R}_\mathbf{d}(\La)$ in the Zariski topology. In this case, $\mathrm{R}_\mathbf{d}(\mcX)$ is a quasi-affine scheme together with an operation of the reductive group $\mathrm{G}_{\mathbf{d}}$. 
\item[(ii)] \textbf{degeneration-open} if for every dimension vector $\mathbf{d}\in \N_{0}^{Q_0}$, $N \in \overline{\mathrm{G}_\mathbf{d}.M}$ and $N \in \mathrm{R}_\mathbf{d}(\mcX)$ imply $M \in \mathrm{R}_\mathbf{d}(\mcX)$.
\end{itemize} 
\end{dfn}
In \cite{DG}, open subcategories are defined in the context of bocses. 
%The definition of an open subcategory is inspired by \cite{DG}, who defined this in the context of bocses. 

\begin{rem}
\begin{itemize}
    \item[(a)] Observe that (i) implies (ii). 
    \item[(b)] Recall that $\La$ (with quiver $Q$) is representation-finite if and only if there are only finitely many $\mathrm{G}_{\mathbf{d}}$-orbits in $\mathrm{R}_{\mathbf{d}}(\La)$ for every $\mathbf{d}\in \N_0^{Q_0}$ (see \cite[Def. I.4.11]{ASS06}). If $\La$ is representation-finite, open and degeneration-open are equivalent, because every $\mathrm{R}_{\mathbf{d}}(\mcX)$ is the complement of finitely many orbit closures.  
%    But the other direction requires that $\La$-mod is representation-finite in order to obtain the complement of $\mathrm{R}_\mathbf{d}(\mcX)$ in $\mathrm{R}_\mathbf{d}(\La)$ as union of finitely many orbit closures.
    \end{itemize}
\end{rem}
\begin{thm}\label{thm:deg-open-equal-ext-closed}
Let $\La$ be a finite-dimensional basic algebra, and $\mcX \subseteq \La$-mod an additively closed subcategory. The following are equivalent
\begin{itemize}
\item[(1)] $\mcX$ is degeneration-open,
\item[(2)] $\mcX$ is extension-closed. 
\end{itemize}
\end{thm}
\begin{proof}
 By \cite[Theorem~1]{Zwa00}, $N \in \overline{\mathrm{G}_\mathbf{d}.M}$ is equivalent to the existence of a short exact sequence $N  \rightarrowtail M \oplus V \twoheadrightarrow V$ for some $V \in \La$-mod. 
 We assume $N$ in $\mcE$. 
In the proof of \cite[Theorem~1]{Zwa00}, Zwara constructs a family of modules $N_i$, $i \geq 1$ , $N_1=N$ with short exact sequences $N_1 \rightarrowtail N_{i+1} \twoheadrightarrow N_i$ for all $i \geq 1$ such that $V\cong N_i$ for some $i>>1$. Now, using these short exact sequences for $i=1$, $N_1\rightarrowtail N_2\twoheadrightarrow N_1$ we conclude that $N_2$ lies in $\mcE$. For $i>1$, inductively, if $N_1,N_i$ are in $\mcE$, then also $N_{i+1}$ is in $\mcE$. Therefore, $V$ can be chosen as an object in $\mcE$. But this implies that $V\oplus M$ lies in $\mcE$. As $\mcE$ is additively closed, it follows that $M$ lies in $\mcE$. Hence, extension-closed implies degeneration-open. For the other direction consider a short exact sequence $U\rightarrowtail V \twoheadrightarrow W$ with $U,W \in \mcE$. This implies $U\oplus W \in \overline{\mathrm{G}_\mathbf{d}.V}$ and hence $V\in \mcE$ if $\mcE$ is degeneration-open.
\end{proof}
\begin{cor} \label{openImpliesExt-closed}
Let $\La$ be a finite-dimensional basic algebra, and $\mcX \subseteq \La$-mod an additively closed subcategory. 
If $\mcX$ is open then it is extension-closed.
\end{cor}
\begin{cor}
    If $\La$ is representation-finite, and $\mcX\subseteq \La$-mod a full additively closed subcategory, the following are equivalent
    \begin{itemize}
    \item[(1)] $\mcX$ is open,
    \item[(2)] $\mcX$ is extension-closed.
    \end{itemize}
\end{cor}

\begin{dfn}
Let $X$ be a $\mathbb{K}$-variety, i.e. the $\mathbb{K}$-valued points of a reduced, separated scheme of finite type over ${\rm Spec} (\mathbb{K})$. A map $f \colon X \to \mathbb{Z}$ is \textbf{upper semi-continuous} if the subset $\{x\in X\, : \, f(x) \leq n\}$ is open (in the Zariski topology) for every $n\in \mathbb{Z}$. A map $f$ is \textbf{lower semi-continuous} if and only if $-f$ is upper semi-continuous.
\end{dfn}

\begin{lem}\label{ext-usc}%(proof or reference...)
Let $\La$ be a finite-dimensional basic algebra with quiver $Q$, $X,Y\in\La$-mod, $\mathbf{d} \in \mathbb{N}_0^{Q_0}$ and $i \in \mathbb{Z}$. %We write $\Ext^0_{\La}:=\Hom_{\La}$. 
Then the map 
\[
\begin{aligned}
\left[X, - \right]^i_{\La}\colon \mathrm{R}_\mathbf{d}(\La) &\to \mathbb{N}_0\\
Y & \mapsto \left[ X, Y \right]^i_{\La}:=\dim_{\mathbb{K}}\Ext^i_{\La} (X,Y)
\end{aligned}
\]
is upper semi-continuous. The same is true for $\left[ -, Y \right]^i_{\La}\colon X\mapsto \dim_{\mathbb{K}}\Ext^i_{\La} (X,Y)$. 
\end{lem}
We include a proof of this statement, as we could not find a suitable reference.
\begin{proof}
By definition (see \cite[p.426]{ASS06}), $\Ext^i_{\La} (X,Y) = \mathrm{Ker} \ h_{i+1}^*/\mathrm{Im} \ h_i^*$ where
\[ \dots \longrightarrow P_{i+1} \overset{h_{i+1}}{\longrightarrow}P_i \overset{h_i}{\longrightarrow} P_{i-1} \longrightarrow \dots\longrightarrow P_1 \overset{h_1}{\longrightarrow} P_0 \overset{h_0}{\longrightarrow} X \longrightarrow 0\]
is a projective resolution of $X$, and $h_{i+1}^* : \Hom_\La(P_{i},Y) \to \Hom_\La(P_{i+1},Y)$ denote the maps in the induced cochain complex. Hence, $\left[ X, Y \right]^i_{\La} = \dim_\K \Hom_\La(P_{i},Y) - \mathrm{rank} \ h_{i+1}^* - \mathrm{rank} \ h_i^*$ is upper semi-continuous, because rank functions on matrices are lower semi-continuous, and each $h_i^*$ is a linear map between vector spaces. For $\left[ -, Y \right]^i_{\La}$, the same arguments apply to the definition of extensions via injective resolutions of $Y$.
\end{proof}

Lemma \ref{ext-usc} provides many examples of open subcategories. 
\begin{exa} \label{ex.open-subcat}
\begin{itemize}
\item[(a)] Fix an object $M$ of $\La$-mod. The subcategory 
\[
{}^{\perp} M:=\{ X\in \La\text{-mod } \, : \, \left[X,M\right]^0_{\La}=0\}
\]
is open and dually $M^{\perp}=\{X\, : \, \left[ M,X\right]^0_{\La}=0\}$ is open. %\\
Observe that finite intersections of open subcategories are open.
Therefore, for any finite subset $I \subseteq \mathbb{N}_0$
\[
{}^{\perp_I}M :=\{ X\in \La\text{-mod } \, : \, \left[X,M\right]^i_{\La}=0 \; \forall i \in I\}
\]
is open (and resp. $M^{\perp_I}$ is open, too). 
\item[(b)] For every $n \in \mathbb{N}_0$, 
the subcategory $\mcP^{\leq n} (\La):=\{X\in \La\text{-mod }  \, : \, \pdim_{\La} X\leq n \}$ is open. To see this: A $\La$-module $X$ has projective dimension $\leq n$ if and only if $\left[X,S\right]^{n+1}_{\La}=0$ for all simple $\La$-modules $S$. Let $S_1, \ldots , S_r$ be all simple $\La$-modules. %(up to isomorphism). 
Then $\mcP^{\leq n} (\La) = \bigcap_{i=1}^r \{ X \, : \, \left[X, S_i\right]^{n+1}_{\La}=0\} $ is a finite intersection of opens, therefore open. %\\
Analogously, $\mcI^{\leq n} (\La ):=\{X \in \La\text{-mod } \, : \, \idim_{\La} X \leq n\}$ is open. 
\end{itemize}
\end{exa}

Furthermore, as a corollary of \cite[Thm 1.3 (iii)]{CS02}, we have the following construction rule: For a class of objects $\mcX$ in $\La\text{-mod}$ we denote by ${\rm Filt} (\mcX)$ the full subcategory of $\La\text{-mod}$ given by modules $X$ which possess a sequence of submodules 
\[
0=X_0\subseteq X_1 \subseteq \cdots \subseteq X_r =X
\]
for some $r \in \mathbb{N}_0$ such that $X_i/X_{i-1}$ is in $\mcX$ for $1\leq i \leq r$. For any full subcategory $\mcY$ of $\La\text{-mod}$, we denote by $\add (\mcY)$ its closure under direct summands of finite direct sums. 

\begin{pro} \label{filt}
Assume that $\mcX_1, \mcX_2, \ldots , \mcX_r$ are open subcategories and $\Ext^1_{\La}(X_j,X_i)=0$ for all $(X_i, X_j)$ in $\mcX_i\times \mcX_j$ for all $1 \leq i < j\leq r$. Then $\add({\rm Filt} (\mcX_1, \ldots , \mcX_r))$ is also open. 
\end{pro}

\begin{proof} Since ${\rm Filt} (\mcX_1, \ldots , \mcX_r) = {\rm Filt} (\mcX_r, {\rm Filt}(\mcX_2, \ldots , \mcX_{r-1}))$, by induction it is enough to show the case $r=2$. Following \cite{CS02} we define for a subset $S\subseteq \mathrm{R}_{\mathbf{d}}(\La)\times \mathrm{R}_{\mathbf{e}}(\La)$ 
\[\mcE (S) := \big\{M \in \mathrm{R}_{\mathbf{d}+\mathbf{e}} (\La)\, : \, \exists \text{ short exact sequence }  M_2 \rightarrowtail M \twoheadrightarrow M_1 \text{ with } (M_1,M_2) \in S \big\}\]
%\[\mcE (S) := \big\{M \in \mathrm{R}_{\mathbf{d}+\mathbf{e}} (\La)\, : \, \exists \text{ ses }  M_2 \rightarrowtail M \twoheadrightarrow M_1 \text{ with } (M_1,M_2) \in S \big\}\]
By \cite[Thm 1.3 (iii)]{CS02}, $\mcE(S)$ is open if $\Ext^1_{\La} (M_2, M_1)=0$ for all $(M_1, M_2) \in S$, and if $S=S_1\times S_2$ we also have $\mcE (S_2\times S_1)= \mathrm{R}_{\mathbf{d}+\mathbf{e}} (S_2 \oplus S_1) \subset \mcE (S)$. %\\

Now, let $\mcX_1, \mcX_2$ be open and $\mcX= {\rm Filt} (\mcX_1, \mcX_2)$ with $\Ext^1_{\La} (X_2, X_1) =0$ for all $(X_1, X_2) \in {\rm Ob}(\mcX_1)\times {\rm Ob}( \mcX_2)$. %We write $\mathbb{d}' \leq \mathbf{d}$ if $\mathbf{d}-\mathbf{d}' \in \mathbb{N}_0^{Q_0}$. 
Then
\[
\mathrm{R}_{\mathbf{d}} (\mcX) = \bigcup_{\mathbf{d}= \mathbf{e}+\mathbf{f}} \mcE (\mathrm{R}_{\mathbf{e}} (\mcX_1) \times \mathrm{R}_{\mathbf{f}}(\mcX_2))
\]
is a finite union of open subsets (using that $\mcE (\mathrm{R}_{\mathbf{f}}(\mcX_2)\times \mathrm{R}_{\mathbf{e}} (\mcX_1))\subset \mcE (\mathrm{R}_{\mathbf{e}} (\mcX_1) \times \mathrm{R}_{\mathbf{f}}(\mcX_2))$ by the Ext-vanishing assumption). Therefore, $\mcX$ is open. The rest of the claim follows from the next lemma. 
\end{proof}

\begin{lem} Let $\mcY$ be a full additive subcategory of $\Lamod$ such that $\mathrm{R}_{\mathbf{d}}(\mcY)$ is open for all dimension vectors $\mathbf{d} \in \mathbb{N}_0^{Q_0}$. Then $\mcX:=\add(\mcY)$ is open. 
\end{lem}

\begin{proof}
We consider the closed immersion $i_{\mathbf{d}\oplus \mathbf{e}} \colon \mathrm{R}_{\mathbf{d}} (\La) \times \mathrm{R}_{\mathbf{e}} (\La) \to \mathrm{R}_{\mathbf{d}+\mathbf{e}} (\La)$ given by $(X,Y)\mapsto X\oplus Y$ and the projection %onto the first factor 
$p_1\colon \mathrm{R}_{\mathbf{d}} (\La) \times \mathrm{R}_{\mathbf{e}} (\La) \to \mathrm{R}_{\mathbf{d}} (\La)$, $(X,Y)\mapsto X$. 
Hence, 
$U_{\mathbf{d},\mathbf{e}}:=p_1(i_{\mathbf{d}\oplus \mathbf{e}}^{-1}(\mathrm{R}_{\mathbf{d}+\mathbf{e}}(\mcY)))$ is an open subset in $\mathrm{R}_{\mathbf{d}}(\La)$, as $p_1$ is an open map. By definition 
\[
\mathrm{R}_{\mathbf{d}}(\mcX) = \bigcup_{\mathbf{e}} U_{\mathbf{d},\mathbf{e}}
\]
is a union of open subsets and therefore is also open. 
\end{proof}

\begin{exa} 
\begin{itemize}
\item[(a)] Given $M_1, \ldots, M_r$ in $\Lamod$ with $\Ext^1(M_j, M_i)=0$ for all $1\leq i \leq j \leq r$, then $\add({\rm Filt} (M_1, \ldots , M_r))$ is open (cf. \cite[Corollary~1.5]{CS02}). For example, $\add({\rm Filt} (I,P))$ is open for every projective $P$ and injective $I$. %\\

An explicit construction of such a collection is found as follows: 
Fix a numbering $S_1, \ldots , S_n$ of the simple $\La$-modules, let $P_1, \ldots , P_n$ be their respective projective cover. The $i$-th standard module is the maximal quotient of $P_i$ with composition factors $S_j$, $j \leq i$. Modules filtered by standard modules is an open subcategory (this does not require any further assumptions such as quasi-heredity), cf.  e.g. \cite{XI_2002}.
\item[(b)] Let $M_1, \ldots ,M_r$ in $\Lamod$ as in (a) with $\Ext^1(M_j, M_i)=0$ for all $1\leq i \leq j \leq r$ and let $M=\bigoplus_{1\leq i\leq r} M_i$. Then we have 
\[
\add ({\rm Filt} (M_1, \ldots , M_r, M^{\perp_1})) \quad \text{and} \quad \add ({\rm Filt} ({}^{\perp_1}M, M_1, \ldots , M_r)) 
\]
are open. For example, for every open subcategory $\mcX$ and every projective $P$ and injective $I$ we have 
$\add ({\rm Filt} (P, I,\mcX))$ is also open.  
\end{itemize}
\end{exa}

%\begin{exa}\label{ex.open-subcat}
%    Let $\mcE \subseteq \La$-mod be a full subcategory described by Hom- or Ext-vanishing conditions. 
%    Then $\mcE \subseteq \La$-mod is open, because $\dim \mathrm{Hom}_\La(-,-)$ and $\dim \mathrm{Ext}^1_\La(-,-)$ are upper semi-continuous \cite[Lemma~4.3]{CS02}. The dimension function of higher extensions is computed via certain rank functions on matrices. Hence, it is also upper semi-continuous and its vanishing can be used to model open subcategories. Analogous to \cite[Corollary~1.5]{CS02}, openness extends to filtered subcategories of $\La$-mod if each quotient of the filtration satisfies some Hom- or Ext-vanishing condition and additionally there are no extensions between a quotient and any of its succeeding quotients.  
%\end{exa}
\section{The Grassmannian of an open fully exact subcategory}
\label{sec:exact-quiver-grass}
From now on assume that $\mcE$ is an open fully exact subcategory of $\Lamod$.
\begin{dfn}
For $\mathbf{e},\mathbf{d} \in \N_{0}^{Q_0}$ and $M \in \mathrm{R}_{\mathbf{e}+\mathbf{d}}(\mathcal{E})$, the $\mathcal{E}$-\textbf{Grassmannian} $\mathrm{Gr}_\mathcal{E}(\mathbf{e},M)$ is the set
\[ \Big\{ U = (U_i)_{i \in Q_0} \in \prod_{i\in Q_0}{\mathrm Gr}_{\K}(e_i,M_i) \ : \ M\vert_U \in \mathrm{R}_{\mathbf{e}}(\mathcal{E}) \ \mathrm{and} \ M/U \in \mathrm{R}_{\mathbf{d}}(\mathcal{E}) \Big\}. \]
\end{dfn}
Our next goal is to equip this set-theoretic Grassmannian with the structure of an algebraic variety. We follow the approach of \cite[Section~4]{CaRe08}. Let $\mathrm{SES}_\mathcal{E}(\mathbf{e},\mathbf{d})$ denote the subset of \[
\mathrm{R}_{\mathbf{e}}(\mathcal{E})\times\prod_{i \in Q_0} \K^{(e_i+d_i)\times e_i} \times \mathrm{R}_{\mathbf{e}+\mathbf{d}}(\mathcal{E})\times \prod_{i \in Q_0} \K^{d_i\times (e_i+d_i)} \times \mathrm{R}_{\mathbf{d}}(\mathcal{E})\] consisting of tuples $(U,f,M,g,N)$ such that each $f_i : \K^{e_i} \to \K^{e_i+d_i}$ is injective, each $g_i : \K^{e_i+d_i} \to \K^{d_i}$ is surjective, $\mathrm{Ker} \, g_i = \mathrm{Im} \, f_i$, and $M_\alpha f_i =  f_j U_\alpha $ and $N_\alpha g_i = f_j M_\alpha $ for all $(\alpha : i \to j) \in Q_1$. The variety structure on $\mathrm{SES}_\mathcal{E}(\mathbf{e},\mathbf{d})$ is obtained from the variety structure on $\mathrm{R}_{\mathbf{d}}(\mathcal{E})$ in combination with the equations arising from the above conditions. For a similar construction, see \cite[Section~3]{CS02}. Here, it is crucial to work with $5$-tuples instead of the middle triples to ensure that all three representations in the short exact sequences are in $\mathcal{E}$ and not only in $\Lambda$-mod. Consider the projections
\begin{center}
 \begin{tikzpicture}[scale=.8]
    \node at (0,2) {$(U,f,M,g,N)$};
    \node[rotate=-90] at (0,1.5) {$\in$};
    \node at (0,1) {$\mathrm{SES}_\mathcal{E}(\mathbf{e},\mathbf{d})$};
    \node at (0,-1) {$ \mathrm{R}_{\mathbf{e}+\mathbf{d}}(\mathcal{E})$};
    
    \node at (-6,-1) {$(U,f,M)$};
    \node at (6,-1) {$(M,g,N)$};

    \node at (-5,-2) {$\mathrm{R}_{\mathbf{e}}(\mathcal{E})\times\prod_{i \in Q_0} \K^{(e_i+d_i)\times e_i} \times \mathrm{R}_{\mathbf{e}+\mathbf{d}}(\mathcal{E})$};
    \node at (5,-2) {$\mathrm{R}_{\mathbf{e}+\mathbf{d}}(\mathcal{E})\times \prod_{i \in Q_0} \K^{d_i\times (e_i+d_i)} \times \mathrm{R}_{\mathbf{d}}(\mathcal{E})$};
    %\node at (-3.164,0) {$M:=$};

    \draw[arrows={-angle 90},shorten >=8, shorten <=8] (0,1) -- (0,-1);

    \draw[arrows={-angle 90},shorten >=25, shorten <=18] (0,1) -- (-5,-2);
    \draw[arrows={-angle 90},shorten >=20, shorten <=18,|->] (-1,2) -- (-6,-1);
    
    \draw[arrows={-angle 90},shorten >=15, shorten <=18] (0,1) -- (5,-2);
    \draw[arrows={-angle 90},shorten >=20, shorten <=18,|->] (1,2) -- (6,-1);

    \node[rotate=30] at (-2.9,-0.25) {$\pi_{infl}$};
    
    \node[rotate=-30] at (2.9,-0.25) {$\pi_{defl}$};
    
    \node at (0.72,-0.35) {$\pi_{\mathbf{e}+\mathbf{d}}$};
    
\end{tikzpicture}
\end{center}
\begin{lem} The maps $\pi_{infl}$ and $\pi_{defl}$ are open. 
\end{lem}
We add an index $\mcE$ to emphasize the dependence on the exact category and an index $\La$ if $\mcE=\Lamod$.

\begin{proof} For $\mcE=\Lamod$, one can stratify $\mathrm{SES}_\mathcal{E}(\mathbf{e},\mathbf{d})$  into locally closed subsets, by fixing $\dim_{\mathbb{K}} \Hom_{\La} (M,N)$. Hence, the restricted map is open by \cite[Lem.2.1]{Bo96}. This implies the claim for $\mcE=\Lamod$. For an open subcategory $\mcE\subseteq \Lamod$, $\mathrm{SES}_\mathcal{E}(\mathbf{e},\mathbf{d})$ is open in $\mathrm{SES}_{\La}(\mathbf{e},\mathbf{d})$ and therefore the claim follows. 
\end{proof}

This implies that the image of $\pi_{infl_{\mcE}}$, denoted by $\mathrm{Infl}_\mathcal{E}(\mathbf{e},\mathbf{d}+\mathbf{e})$ is an open subset of $\mathrm{Infl}_{\La }(\mathbf{e},\mathbf{d}+\mathbf{e})$. In particular  $\mathrm{Infl}_\mathcal{E}(\mathbf{e},M) = \pi_{infl_{\mcE}}(\pi^{-1}_{\mathbf{e}+\mathbf{d}}(M))$ is open and $\mathrm{G}_\mathbf{e}$-invariant in $\mathrm{Infl}_{\La}(\mathbf{e},M) = \pi_{infl_{\La}}(\pi^{-1}_{\mathbf{e}+\mathbf{d}}(M))$. 
% $\mathrm{Infl}_\mathcal{E}(\mathbf{e},M) = \pi_{infl}(\pi^{-1}_{\mathbf{e}+\mathbf{d}}(M))$ and $\mathrm{Defl}_\mathcal{E}(\mathbf{e},M) = \pi_{defl}(\pi^{-1}_{\mathbf{e}+\mathbf{d}}(M))$ as it contains the $\mathbf{e}$-dimensional inflations (resp. $\mathbf{d}$-dimensional deflations) of $M$. 
%The maps in the analogous construction for $\La$-mod are fiber bundles and this is preserved passing to $\mcE$-representations. Hence, the variety structure on  $\mathrm{Infl}_\mathcal{E}(\mathbf{e},M)$ and $\mathrm{Defl}_\mathcal{E}(\mathbf{e},M)$ is obtained from $\mathrm{SES}_\mathcal{E}(\mathbf{e},\mathbf{d})$. 
Here, $\mathrm{G}_\mathbf{e}$ acts freely on $\mathrm{Infl}_\mathcal{E}(\mathbf{e},M)$ via 
\[g.(U,f) := \big( (g_jU_\alpha g_i^{-1})_{(\alpha:i \to j) \in Q_1}, (f_i g_i^{-1})_{i \in Q_0} \big).\]
The following result generalizes \cite[Lemma~2]{CaRe08}. %The proof is based on ideas from \cite{CaRe08} and \cite[Section~2.3]{CFR11}. 
For more detail also see %\cite[Proposition~2.5]{Pu19} and 
\cite[Section~3.2]{Re08}.

\begin{lem}\label{lem:subvariety-structure} 
Let $\mcE$ be an open fully exact subcategory of $\Lamod$. 
    The $\mathcal{E}$-Grassmannian $\mathrm{Gr}_\mathcal{E}(\mathbf{e},M)$ is isomorphic to the geometric quotient
    $\mathrm{Infl}_\mathcal{E}(\mathbf{e},M)/\mathrm{G}_\mathbf{e}$. In particular, it is an open subvariety of $\mathrm{Gr}_\Lambda(\mathbf{e},M)$.
\end{lem}
\begin{proof} 
    Since $\phi_{\La} : \mathrm{Infl}_{\La}(\mathbf{e},M) \to  \mathrm{Gr}_\La(\mathbf{e},M)$ sending $(U,f)$ to $\oplus_{i\in Q_0} f_i(\K^{e_i})$ is a geometric quotient, the same holds for its restriction to an open $\mathrm{G}_\mathbf{e}$-invariant subscheme. 
%    $\phi : \mathrm{Infl}_\mathcal{E}(\mathbf{e},M) \to  \mathrm{Gr}_\mathcal{E}(\mathbf{e},M)$ sending $(U,f)$ to $\oplus_{i\in Q_0} f_i(\K^{e_i})$. The image is a submodule of $M$ since $M_\alpha f_i = f_j U_\alpha$ holds. Hence $\phi$ is well defined and clearly it is surjective. The injectivity of the $f_i$'s implies that $\phi$ is constant on $\mathrm{G}_\mathbf{e}$-orbits in $\mathrm{Infl}_\mathcal{E}(\mathbf{e},M)$. Now, the variety structure on the quotient is uniquely determined by $\mathcal{O}_{\mathrm{Gr}_\mathcal{E}(\mathbf{e},M)} = (\phi_*\mathcal{O}_{\mathrm{Infl}_\mathcal{E}(\mathbf{e},M)})^{\mathrm{G}_\mathbf{e}}$, where it is proved as in \cite[Lemma~3.5]{Re08} that the quotient map is algebraic. %\cite[Section~2.3]{CFR11}. 
\end{proof}
\begin{rem}
    Based on deflations, there is a dual construction of 
    the $\mathcal{E}$-Grassmannian $\mathrm{Gr}_\mathcal{E}(M,\mathbf{d})=\mathrm{Defl}_\mathcal{E}(M,\mathbf{d})/\mathrm{G}_\mathbf{d}$. In particular, $\mathrm{Gr}_\mathcal{E}(\mathbf{e},M) \cong \mathrm{Gr}_\mathcal{E}(M,\mathbf{d})$ via $U \mapsto M/U$. 
\end{rem}
%++++++++++++++++++++++++++++++++++++++++++++++++++++++++++++
\section{Quasi-projective varieties and Grassmannians for exact categories}\label{sec:quasi-proj}
%++++++++++++++++++++++++++++++++++++++++++++++++++++++++++++
We start by recalling Reineke's isomorphism (from \cite{R}). % and first look at a single principle open subset of a projective variety. 
 %\noindent
%Let $\K$ be an algebraically closed field. 
We consider $X= \rm{Proj} (R) $ for $R=\K[T_0, \ldots , T_n]/(f_1, \ldots , f_k)$, with homogeneous polynomials $f_i$ of degree $d(f_i)>0$, for $1 \leq i \leq k$.  Let $h_1, \ldots , h_s$ be further homogeneous polynomials such that $\overline{h_i} \in R$ has a positive degree $d(h_j)>0$. 
In particular, $h_j$ is not a scalar multiple of any $f_i$, for $1 \leq i \leq k$. Then the open subset $\bigcup_{j=1}^s D_+(\overline{h_j})\subset X$ equals 
\[ \big\{x\in \mathbb{P}^n(\K) \, : \, f_i(x)=0 \text{ for all } 1 \leq i \leq k\text{ and there exists }  1 \leq j \leq s \colon  h_j(x)\neq 0\big\}.\]
%\[ \{x\in \mathbb{P}^n(\K) \, : \, \forall \; i \in [k] \colon  f_i(x)=0\text{ and } \exists j \in [s]\colon  h_j(x)\neq 0\}.\]
    Assume $d:=d(f_i)=d(h_j)$ for all $1 \leq i \leq k$ and all $1 \leq j \leq s$. Otherwise, replace $f_i$ by $f_i^{t_i}$ and $h_j$ by $h_j^{s_j}$ where $\ell = \mathrm{lcm} (d(f_1), \ldots , d(f_k), d(h_1), \ldots , d(h_s))$, $t_i = {\ell}/{d(f_i)}$ and $s_j = {\ell}/{d(h_j)}$. This does not change the variety or the open subset. 
    Furthermore, by allowing repetitions in $f_1, \ldots , f_k$ and $h_1, \ldots , h_s$ respectively, we can assume that $s=k$. 
    
    %\item[(2)] 
    To assign a quiver representation to the equations of the variety $X$, we use the $d$-uple embedding of $\mathbb{P}^n$: This embedding linearizes the describing equations of the variety. Let
     \begin{center}         
     \(M_{n,d}:= \big\{ (m_0, \ldots , m_n) \in \mathbb{N}_0^{n+1} \ : \ \sum_{j=0}^n m_j= d\big\},\)\end{center}
     and let $N:=\lvert M_{n,d} \rvert$ denote its cardinality. Define 
     \[ 
    \begin{aligned}
    \nu \colon \mathbb{P}^n &\to \mathbb{P}^{N-1},\quad [x_0\colon \ldots \colon x_n ] \mapsto [\cdots \colon x_\mathbf{m}\colon \cdots ]_{\mathbf{m} \in M_{n,d}}, \\ &\text{where } \quad x_\mathbf{m} := x_0^{m_0} \cdots x_n^{m_n}, \;\text{for }\mathbf{m} =(m_0, \ldots , m_n) \in M_{n,d}
    \end{aligned}
    \]
    This is a closed embedding called the \textbf{$d$-th Veronese embedding} (or $d$-uple embedding).
    
Now, we consider the quiver $Q=(Q_0, Q_1)$ with three vertices $1,2,3$ and $n+1$ arrows from $2$ to $3$ and $k$ arrows from $2$ to $1$. 
    We define the following $Q$-representations: 
    Let $\{v_\mathbf{m} : \mathbf{m}\in M_{n,d}\}$ be a vector space basis of $\K^{M_{n,d}}$. This is the set of all linear maps $M_{n,d} \to \K$, given a vector space structure with point-wise addition and scalar multiplication. 
    We denote by $\varphi_i$ (resp. $\psi_i $) the linear maps $\K^{M_{n,d}}\to \K$ sending $v_\mathbf{m} \mapsto a_\mathbf{m}^{(i)}$ (resp. $b^{(j)}_\mathbf{m}$), corresponding to $f_i= \sum_{\mathbf{m}\in M_{n,d}}a_\mathbf{m}^{(i)}T^\mathbf{m}$, and $h_i= \sum_{\mathbf{m}\in M_{n,d}} b^{(i)}_\mathbf{m} T^\mathbf{m}$ respectively.
    %\\
    
    Let $V$ be the $Q$-representation with vector spaces $V_1= \K$, $V_2=\K^{M_{n,d}}$, $V_3= \K^{M_{n,d-1}}$ and linear maps $\varphi_i\colon V_2 \to V_1$, for $1 \leq i  \leq k$, and $g_j \colon V_2 \to V_3$, for $0\leq j \leq n$ defined by $g_j(v_\mathbf{m}):= v_{\mathbf{m}-\mathbf{e}_j}$ if $m_j>0$ and $g_j(v_\mathbf{m}):=0$ if $m_j=0$. Here, $\mathbf{e}_j$ denotes the tuple with entry $1$ in $j$-th position. %\\
    We define the $Q$-representation $W$ with the same underlying vector spaces $W_i=V_i$, $i\in \{1,2,3\}$  and $g_j \colon W_2 \to W_3$, $0 \leq j \leq n$ defined as in $V$, but linear maps $\psi_i\colon W_2 \to W_1$, for $1 \leq i \leq k$. %\\ 
    \begin{thm}(Reineke \cite{R})
    The map $\Gr_{Q} (\mathbf{d},V) \to \nu(X)\cong X, (i \colon U \rightarrowtail V) \mapsto \Bild i_2$ is an isomorphism, for $\mathbf{d}=(0,1,1)$.
    \end{thm}
    %In \cite{R} it is shown that $\nu(X)$ is isomorphic to $\Gr_Q(\mathbf{d},V)$, for $\mathbf{d}=(0,1,1)$. %\\
    
    For every $x=(x_0, \ldots , x_n)\in \mathbb{K}^{n+1}\setminus{\{0\}}$ we define a representation $U_x$ of $Q$:
    \[
    \xymatrix{ 
    0 & \mathbb{K}\ar@<2ex>[r]^{x_0}\ar@<-1ex>[r]_{x_n}^\vdots  \ar@<1ex>[l]_{\vdots} \ar@<-2ex>[l] & \mathbb{K} 
    }
    \]
    Two monomorphisms $U_x\to V$, $U_{x'}\to V$ have the same image (i.e. are the same point in the Grassmannian) if and only if $[x]=[x']\in \mathbb{P}^n$. We write $U_{[x]}$ for the image (i.e. the submodule of $V$), which is as $Q$-representation isomorphic to $U_x$. 
    
    \begin{lem}\label{lem:some-hom-results}
    \begin{itemize}
    \item[(i)] $\Hom_{Q} (V, W)=0= \Hom_{Q} (W, V)$ and $\End_{Q}(W)=\K= \End_{Q} (V)$.
    \item[(ii)] If $U_{[x]} \in \Gr_{Q} (\mathbf{d},V)$ such that $i\colon U_x \rightarrowtail V$ corresponds
    to $[x]\in X\subset \mathbb{P}^n(\K)$ (i.e. $\nu([x])= \Bild i_2$), then: %\\
    $\Hom_{Q} (U_x, W)=0$ if and only if there exists $1\leq j \leq k$ such that $h_j(x)\neq 0$. 
    \end{itemize}
    \end{lem}
    
    \begin{proof}
        (i): Let $Q'$ be the full subquiver on the vertices $\{2,3\}$. The restricted representations onto this subquiver coincide, i.e. $V':=V\vert_{Q'}=W\vert_{Q'}$. 
        In the last paragraph of \cite{R} it is shown that $\End_{Q'} (V')=\K$. More precisely, every endomorphism consists of a pair of linear maps $\Psi_2 \colon V_2 \to V_2$, $\Psi_3\colon V_3 \to V_3$ with $\Psi_2 = C\id_{V_2}$ for $C \in \K$ and $\Psi_3$ is determined by $\Psi_2$. Now, since no $h_j$ is a scalar multiple of the $f_i$'s, the rest of the claim is clear from the definitions.
        
        (ii): Following the last paragraph of \cite{R},  $\Hom_{Q} (U_x, W)\neq 0$ is equivalent to the existence of a monomorphism $U_x \rightarrowtail W$. By construction of $U_x$ and $W$, this is equivalent to $h_j(x)=0$ for all $j$.%, as $i\colon U_x \rightarrowtail V$ exists by assumption.
    \end{proof}
    
    Let $\mcE ={}^\perp W\subset \K Q\mathrm{-mod}$ be the full subcategory consisting of all modules $X$ with $\Hom_{Q} (X, W)=0$. This is an open subcategory, cf. Example \ref{ex.open-subcat}, and we consider it equipped with the fully exact structure. Given a monomorphism $U_x\rightarrowtail V$, then $V\in \mcE$ and $V/U_x\in \mcE$ are always fulfilled (by Lemma \ref{lem:some-hom-results} (i) and as $\mcE$ is closed under factor modules). Therefore, $U_x\rightarrowtail V \twoheadrightarrow V/U_x$ is a short exact sequence in $\mcE$ if and only if $U_x \in \mcE$, i.e. $\Hom_{Q} (U_x, W)=0$. The previous lemma directly implies:

\begin{thm}\label{thm:quasi-projective-varieties-as-quiver-Grass}
  Reineke's isomorphism $\Gr_{Q} (\mathbf{d},V) \to \nu(X)\cong X, (i \colon U \rightarrowtail V) \mapsto \Bild i_2$, restricts to an isomorphism of open 
    subvarieties $\Gr_{\mcE} (\mathbf{d},V)\to \bigcup_{j =1}^k D_+(\overline{h}_j)$ where $\mcE$ is the open fully exact category ${}^\perp W$ as before. 
    \end{thm}
As our field is algebraically closed, every brick is a Schurian representation and vice versa (i.e. a representation with endomorphism ring isomorphic to $\mathbb{K}$). Then $V,W$ are Hom-orthogonal bricks by Lemma \ref{lem:some-hom-results}. 

\begin{rem}
    The map is algebraic for the same reasons as in Reineke's setting. 
\end{rem}

%\bibliographystyle{amsalpha}
%\bibliography{Lit}
  
\end{document}